\documentclass[11pt]{article}
\usepackage{amsmath}
\usepackage{amsfonts}
\usepackage{color}
\usepackage[utf8]{inputenc}
\usepackage[T1]{fontenc}
\usepackage{hyperref}
\hypersetup{
colorlinks=true,
linkcolor=blue,
citecolor=blue,
urlcolor=blue,
pdftitle={Chemotaxis-consumption survey},
pdfauthor={Lankeit, Winkler},
}
\newtheorem{theo}{Theorem}[section]

\newtheorem{prop}[theo]{Proposition}

\newcommand{\mysection}[1]{\section{#1} \setcounter{equation}{0}}
\newcommand{\proof}{{\sc Proof.} \quad}

\newcommand{\be}{\begin{equation} \label}
\newcommand{\ee}{\end{equation}}
\newcommand{\bea}{\begin{eqnarray}\label}
\newcommand{\eea}{\end{eqnarray}}
\newcommand{\bas}{\begin{eqnarray*}}
\newcommand{\eas}{\end{eqnarray*}}
\newcommand{\bit}{\begin{itemize}}
\newcommand{\eit}{\end{itemize}}
\newcommand{\qed}{\hfill$\Box$ \vskip.2cm}
\newcommand{\nn}{\nonumber}
\newcommand{\R}{\mathbb{R}}

\newcommand{\pO}{\partial\Omega}

\newcommand{\supp}{{\rm supp} \, }

\newcommand{\io}{\int_\Omega}
\newcommand{\na}{\nabla}
\newcommand{\Del}{\Delta}
\newcommand{\del}{\delta}
\newcommand{\al}{\alpha}
\newcommand{\lam}{\lambda}
\newcommand{\sig}{\sigma}
\newcommand{\pa}{\partial}
\newcommand{\bom}{\overline{\Omega}}
\newcommand{\Om}{\Omega}

\newcommand{\wh}{\widehat}

\newcommand{\vs}{\vspace*}
\newcommand{\hs}{\hspace*}

\newcommand{\vp}{\varphi}
\newcommand{\lbal}{\left\{ \begin{array}{l}}
\newcommand{\lball}{\left\{ \begin{array}{ll}}
\newcommand{\ear}{\end{array} \right.}
\newcommand{\abs}{\\[5pt]}

\newcommand{\F}{{\mathcal{F}}}
\newcommand{\D}{{\mathcal{D}}}

\newcommand{\ouz}{\overline{u}_0}
\usepackage{newunicodechar}
\newunicodechar{μ}{\mu}
\newunicodechar{Δ}{\Delta}
\newunicodechar{ℝ}{\R}
\newunicodechar{δ}{\delta}
\newunicodechar{β}{\beta}
\newunicodechar{ν}{\nu}
\begin{document}
\enlargethispage{10mm}
\author{Johannes Lankeit\footnote{lankeit@ifam.uni-hannover.de}\\
\small{Leibniz Universität Hannover}\\\small{ Institut für Angewandte
  Mathematik}\\\small{ Welfengarten 1}\\\small{ 30167 Hannover, Germany}
\and
Michael Winkler\footnote{michael.winkler@math.uni-paderborn.de}\\
{\small Institut f\"ur Mathematik}\\{\small Universit\"at Paderborn}\\
{\small Warburger Str. 100}\\\small{ 33098 Paderborn, Germany} 
}  

\title{Depleting the signal: Analysis of chemotaxis-consumption models -- A survey}

\date{}
\maketitle
\begin{abstract}
\noindent 
We give an overview of analytical results concerned with chemotaxis systems where the signal is absorbed. 
We recall results on existence and properties of solutions for the prototypical chemotaxis-consumption model and various variants 
and review more recent findings on its ability to support the emergence of spatial structures.
\abs
\end{abstract}
\newpage
\mysection{Introduction}\label{intro}
Many microbes are motile, 
and often the direction of their motion may be influenced by an external cue. If said stimulus is provided by the concentration of a chemical signal substance, e.g., serving as a nutrient, this phenomenon is known as \textit{chemotaxis}. 

Already in 1881, Engelmann used bacteria chemotactically attracted by oxygen to visualize production of the latter \cite{engelmann}; more detailed investigations with focus on chemotaxis itself were reported by Pfeffer \cite{pfeffer1,pfeffer2}, and Beyerinck in 1893 demonstrated macroscopically visible manifestations of this effect  
\cite{beyerinck}. 

Experiments by Adler (in the 1960's, \cite{adler66,adler_effect_of}) to quantify and to better understand the mechanism of chemotaxis in bacteria motivated Keller and Segel \cite{keller_segel_traveling} to ``formulate a phenomenological model from which the existence and properties of migrating bands can be deduced'' \cite[p.236]{keller_segel_traveling}.

In an even simpler form, with prototypical choices for all parameter functions and all constants set to $1$, a PDE system for the evolution of the distribution of bacteria (with density $u$) and the signal (of concentration $v$) reads 
\be{CC1}
	\lball
	u_t=\Del u - \na\cdot (u\na v), 
	\qquad & x\in\Om, \ t>0, \\[1mm]
	v_t = \Del v - uv
	\qquad & x\in\Om, \ t>0, \\[1mm]
	\frac{\pa u}{\pa\nu}=\frac{\pa v}{\pa\nu}=0,
	\qquad & x\in\pO, \ t>0, \\[1mm]
	u(x,0)=u_0(x), \quad v(x,0)=v_0(x),
	\qquad & x\in\Om,
	\ear
	\tag{CC1}
\ee
where $\Omega$ usually is a domain in $ℝ^2$ or, especially when considering scenarios inside liquid environments, in $ℝ^3$, and where $ν$ denotes the outer unit normal on its boundary.

In the mathematical literature long overshadowed by its famous cousin, ``the'' Keller--Segel system, 
\be{KS}
	\left\{ \begin{array}{l}
	u_t = \Delta u - \nabla \cdot (u\nabla v),\\[1mm]
	v_t = \Delta v-v+u,
	\end{array} \right.
\ee
where the signal is produced by the cells under consideration and not 'only' consumed, \eqref{CC1} along with its variants received renewed interest, when in the wake of experimental observations concerned with bacteria in drops of water \cite{goldstein2004,goldstein2005}, starting from \cite{lorz10,duan_lorz_markowich} models coupling chemotaxis effects with (Navier-)Stokes fluid motion of their surroundings became popular. 

While compared to \eqref{KS}, \eqref{CC1} does not feature production and thus possible increase of the signal concentration, the boundedness information afforded by the second equation (in \eqref{CC1}, an $L^\infty(\Om)$ bound for $v$ is immediately assured due to the sign of the nonlinearity) is still insufficient to render a study of the first equation trivial. 

In this survey, we aim to collect some of the results that have been achieved for \eqref{CC1} and close relatives over the last years, and to recall some of the underlying ideas. 

\mysection{Classical chemotaxis-consumption systems}
In the apparently most prototypical chemotaxis-consumption problem, as given by \eqref{CC1}, 
the somewhat antagonistic character of the interplay between the two crucial nonlinearities becomes manifest in 
an energy identity of the form
\be{en}
	\frac{d}{dt} \F_1(t) = - \D_1(t),
	\qquad t>0,
\ee
formally associated with (\ref{CC1}).
Here, unlike in a corresponding identity for classical Keller-Segel production systems, not only the dissipation rate
\be{D}
	\D_1(t):=\io \frac{|\na u|^2}{u} + \io v|D^2 \ln v|^2
	- \frac{1}{2} \io \frac{u}{v} |\na v|^2
	- \frac{1}{2} \int_{\pO} \frac{1}{v} \frac{\pa |\na v|^2}{\pa\nu},
	\qquad t>0,
\ee
but also the energy functional (\cite{duan_lorz_markowich,win_CPDE,taowin_consumption})
\be{F}
	\F_1(t):=\io u\ln u + \frac{1}{2} \io \frac{|\na v|^2}{v},
	\qquad t>0,
\ee
can immediately seen to be bounded from below at least when $\Om$ is assumed to be convex, as then 
any function $\vp\in C^1(\bom)$ fulfilling $\frac{\pa\vp}{\pa\nu}=0$ on $\pO$ satisfies $\frac{\pa|\vp|^2}{\pa\nu}\le 0$
throughout $\pO$ (\cite{lions_ARMA}).\abs
Thanks to the functional inequality
\be{FE}
	\io \frac{|\na \vp|^4}{\vp^3} \le (2+\sqrt{n})^2 \io \vp|D^2\ln \vp|^2,
\ee
valid actually for arbitrary bounded domains $\Om\subset\R^n$ with smooth boundary, and for any
positive $\vp\in C^2(\bom)$ such that $\frac{\pa\vp}{\pa\nu}|_{\pO}=0$ (\cite{win_CPDE}), 
in the convex case an integration of (\ref{en}) especially yields estimates of the form
\be{CC1.1}
	\int_0^T \io \frac{|\na u|^2}{u} + \int_0^T \io \frac{|\na v|^4}{v^3} + \int_0^T \io \frac{u}{v} |\na v|^2 \le C
	\qquad \mbox{for all } T>0,
\ee
so that since
\be{vinfty}
	\|v(\cdot,t)\|_{L^\infty(\Om)} \le \|v_0\|_{L^\infty(\Om)}
	\qquad \mbox{for all } t>0
\ee
by the maximum principle, {\em a priori} bounds for $\int_0^\infty \io |\na v|^4$ are available.\abs
In the case when $\Om$ additionally is two-dimensional, this information on regularity of the taxis gradient in (\ref{CC1}) 
can be used as a starting point for a boothstrap procedure finally leading not only to $L^\infty$ bounds for $u$,
but furthermore also to a statement on large time stabilization (\cite{taowin_consumption}).
As observed in \cite{songmu_zheng}, by using boundary trace embedding estimates to appropriately estimate
the rightmost summand in (\ref{D}), an extension to actually arbitrary and not necessarily convex planar domains is possible:
\begin{theo}\label{theo1} (\cite{songmu_zheng})	\quad		
  Let $n=2$ and $\Om\subset\R^n$ be a bounded domain with smooth boundary, and suppose that
  \be{init}
	\lbal
	u_0 \in W^{1,\infty}(\Om) \mbox{ is nonnegative with $u_0\not\equiv 0$, \quad and that} \\[1mm]
	v_0 \in W^{1,\infty}(\Om) \mbox{ is positive in $\bom$.}
	\ear
  \ee
  Then there exist uniquely determined functions
  \be{reg}
	\lbal
	u\in C^0(\bom\times [0,\infty)) \cap C^{2,1}(\bom\times (0,\infty)) 
	\qquad \mbox{and} \\[1mm]
	v\in \bigcap_{q>n} C^0([0,\infty);W^{1,q}(\Om)) \cap C^{2,1}(\bom\times (0,\infty)) 
	\end{array} \right.
  \ee
  such that $u\ge 0$ and $v>0$ in $\bom\times (0,\infty)$, and that (\ref{CC1}) is solved in the classical sense. 
  Moreover, as $t\to\infty$ we have
  \be{conv}
	u(\cdot,t) \to \ouz := \frac{1}{|\Om|} \io u_0
	\quad \mbox{in } L^\infty(\Om)
	\qquad \mbox{and} \qquad
	v(\cdot,t) \to 0
	\quad \mbox{in } L^\infty(\Om).
  \ee
\end{theo}
In its higher-dimensional version, (\ref{CC1}) after all admits some globally defined solutions within weaker
concepts of solvability.
In formulating two corresponding results here, we concentrate on the case of convex domains, noting that extensions
to non-convex situations could be achieved at the cost of some additional technical expense (\cite{songmu_zheng}).
\begin{theo}\label{theo2} (\cite{taowin_consumption}, \cite{li_yuxiang_EJDE}) \quad			
  Let $n\ge 3$ and $\Om\subset\R^n$ be a bounded convex domain with smooth boundary, and suppose that
  $u_0$ and $v_0$ comply with (\ref{init}).\abs
  i) \ If $n=3$, then one can find nonnegative functions
  \be{reg2}
	\lbal
	u\in L^\infty((0,\infty);L^1(\Omega)) 		
	\qquad \mbox{and} \\[1mm]
	v\in L^\infty(\Omega\times (0,\infty)) \cap L^2((0,\infty);W^{1,2}(\Omega))
	\ear
  \ee
  such that 
  \be{2.1}
	\na u \in L^1_{loc}(\bom\times [0,\infty);\R^n)
	\qquad \mbox{and} \qquad
	u\na v\in L^1_{loc}(\bom\times [0,\infty);\R^n),
  \ee
  that $u\ge 0$ and $v>0$ in $\bom\times (0,\infty)$, and that (\ref{CC1}) is solved in the sense that for
  arbitrary $\vp\in C_0^\infty(\bom\times [0,\infty))$ we have
  \be{2.2}
	- \int_0^\infty \io u\vp_t - \io u_0\vp(\cdot,0)
	= - \int_0^\infty \io \na u\cdot\na \vp
	+ \int_0^\infty \io u\na v\cdot\na\vp
  \ee
  and
  \be{2.3}
	\int_0^\infty \io v\vp_t + \io v_0 \vp(\cdot,0)
	= \int_0^\infty \io \na v\cdot\na\vp
	+ \int_0^\infty \io uv\vp.
  \ee
  ii) \ If $n\ge 3$ is arbitrary, then there exist nonnegative functions $u$ and $v$ fulfilling (\ref{reg2}), which are such that
  \bas
	\xi(u)\in L^2_{loc}([0,\infty);W^{1,2}(\Om))
	\qquad \mbox{for all $\xi\in C^\infty([0,\infty))$ with } \xi'\in C_0^\infty([0,\infty)),
  \eas
  and that for any such $\xi$ and arbitrary $\vp\in C_0^\infty(\bom\times [0,\infty))$, both (\ref{2.3}) and the identity
  \bas
	- \int_0^\infty \io -\xi(u)\vp_t - \io \xi(u_0)\vp(\cdot,0)
	&=& - \int_0^\infty \io \xi''(u) |\na u|^2 \vp
	- \int_0^\infty \io \xi'(u) \na u\cdot\na\vp \nn\\
	& & + \int_0^\infty \io u\xi''(u) (\na u\cdot\na v)\vp
	+ \int_0^\infty \io u\xi'(u) \na v\cdot\na\vp
  \eas
  hold.
\end{theo}
To date yet unsolved seems the question whether for $n\ge 3$ and initial data of arbitrary size also global classical
solutions to the fully parabolic system (\ref{CC1}) can be expected 
(cf.~also \ref{theo4} below, \cite{jie_jiang11} for a criterion on extensibility of local-in-time classical solutions,
as well as \cite{taowin_JDE_parab_ell} for
a partial result addressing unconditional global classical solvability in a parabolic-elliptic simplification of (\ref{CC1})). 
On the other hand, the totally inviscid system corresponding to \eqref{CC1}, while locally well-posed, has solutions blowing up in finite time (with respect to their $C^2$-norm), \cite{kyungkeun_CMP}.

Ultimately, the decay features implicitly expressed by (\ref{CC1.1}) imply, at least formally, 
smallness of $\io \frac{|\na u|^2}{u}$ and of $\io |\na v|^4$ along some unbounded sequence of times. 
Now at least in three-dimensional cases, an appropriate exploitation of this can be used to infer that the respective weak solutions
constructed in Theorem \ref{theo2} become smooth and bounded after some waiting time:
\begin{theo}\label{theo3} (\cite{taowin_consumption}) \quad
  Let $\Om\subset\R^3$ be a bounded convex domain with smooth boundary, and assume (\ref{init}).
  Then there exists $T>0$ such that the weak solution $(u,v)$ of (\ref{CC1}) from Theorem \ref{theo2} satisfies
  \be{3.1}
	(u,v) \in \big( C^{2,1}(\bom\times [T,\infty)) \big)^2.
  \ee
  Moreover, this solution stabilizes in the sense that (\ref{conv}) holds.
\end{theo}
Especially since the energy structure in (\ref{en}) naturally is rather fragile with respect to modifications in (\ref{CC1}),
an approach different from the above has turned out to be of significant importance not only for the analysis of (\ref{CC1}),
but also of several more complex relatives. 
As observed in \cite{tao} (cf.~also \cite{win_M2AS2011} for a precedent),	
namely, for suitably chosen $p>1$, $\del>0$ and $\vp:[0,\del] \to (0,\infty)$,
the expressions
\be{F2}
	\F_2(t) := \io u^p \vp(v),
	\qquad t>0,
\ee
can play the role of {\em conditional} quasi-Lyapunov functionals in the sense that with some $a\ge 0, b\ge 0$ and
$\wh{\D}=\wh{\D}(t)\ge 0$,
\be{en2}
	\F_2'(t) + a\F_2(t) + \wh{D}(t) \le b
\ee
holds as long as $0\le v(x,t) \le \del$.\abs
In \cite{tao}, the particular choices $p:=n+1$, $\del:=\frac{1}{6(n+1)}$ and
\be{phi}
	\vp(s):=e^{\beta s^2},
	\quad s\ge 0,
	\qquad
	\mbox{with} \quad \beta:=\frac{n}{24(n+1)\del^2},
\ee
led to the conclusion that whenever (\ref{init}) holds with $\|v_0\|_{L^\infty(\Om)} \le \del$,
then (\ref{en2}) is satisfied throughout evolution if $\wh{D}\equiv 0$, $a=1$ and $b=b(u_0,v_0)>0$ is suitably large.
Accordingly implied time-independent estimates for $u$ with respect to the norm in $L^{n+1}(\Om)$, however, can be seen
to imply higher regularity properties, and to thus finally entail global existence of smooth
small-signal solutions to (\ref{CC1}). The following statement in this direction contains the yet slightly weaker
assumption (\ref{4.1}), obtained by a refinement of the above idea in \cite{baghaei}. (This constant allows for further improvements. For example, one could replace $L^{n+1}(\Om)$ by $L^{\frac{n}{2}+\eta}(\Om)$ for some $\eta>0$, upon observing that time-independent bounds in the latter space already entail the same higher regularity properties as before.)
\begin{theo}\label{theo4} (\cite{tao}, \cite{baghaei})	\quad
  Let $n\ge 2$ and $\Om\subset\R^n$ be a bounded domain with smooth boundary, and suppose that 
  $u_0$ and $v_0$ satisfy (\ref{init}) as well as
  \be{4.1}
	\|v_0\|_{L^\infty(\Om)} \le \frac{\pi}{\sqrt{2(n+1)}}.		
  \ee
  Then (\ref{CC1}) possesses a unique global classical solution fulfilling (\ref{reg}). 
  Moreover, this solution is bounded in the sense that there exists $C>0$ such that
  \be{bound}
	\|u(\cdot,t)\|_{L^\infty(\Om)} \le C
	\qquad \mbox{for all } t>0.
  \ee
\end{theo}
\vs{4mm}
In fact, a further development of this on the basis of an approach from \cite{taowin_consumption} reveals that 
at least for small-signal trajectories even the choice $b=0$ can be achieved in (\ref{en2}) if $\vp$ is properly chosen;
in consequence, any such solution can be seen to approach the semi-trivial steady state appearing in (\ref{conv}):
\begin{prop}\label{prop5} 
  Assume that $n\ge 2$, and that $\Om\subset\R^n$ is a bounded domain with smooth boundary.
  Then there exists $\del\in (0,\frac{\pi}{\sqrt{2(n+1)}}]$ with the property such that if (\ref{init}) holds with 
  \be{5.1}
	\|v_0\|_{L^\infty(\Om)} < \del,
  \ee
  then the global solution of (\ref{CC1}) from Theorem \ref{theo4} satisfies (\ref{conv}).
\end{prop}
\proof
  We let $\del \in (0,\frac{\pi}{\sqrt{2(n+1)}}]$ be such that
  \be{5.2}
	\del\le\frac{1}{8}
	\qquad \mbox{and} \qquad
	32\del + 24\del^2 < \frac{4}{3},
  \ee
  and assuming that (\ref{init}) and (\ref{5.1}), we then take $(u,v)$ as provided by Theorem \ref{theo4} and integrate
  by parts using (\ref{CC1}) to see that since $0\le v < \del$ in $\bom\times [0,\infty)$ by (\ref{vinfty}),
  \bea{5.3}
	\frac{d}{dt} \io \frac{(u-1)^4}{\del-v}
	&=& - 12 \io \frac{(u+1)^2}{\del-v} |\na u|^2
	- \io \frac{(u+1)^3}{(\del-v)^2} \cdot \Big\{ 8 - \frac{12(\del-v)u}{u+1} \Big\} \na u\cdot\na v \nn\\
	& & - \io \frac{(u+1)^4}{(\del-v)^3} \cdot \Big\{ 2 - \frac{4(\del-v)u}{u+1}\Big\} |\na v|^2
	- \io \frac{u^5 v}{(\del-v)^2}
	\qquad \mbox{for all } t>0.
  \eea
  As
  \bas
	2 - \frac{4(\del-v)u}{u+1} \ge 2 - 4\del \ge \frac{3}{2}
	\qquad \mbox{in } \Om\times (0,\infty)
  \eas
  by (\ref{5.2}), 
  we may employ Young's inequality to estimate
  \bas
	& & \hs{-20mm}
	- \io \frac{(u+1)^3}{(\del-v)^2} \cdot \Big\{ 8 - \frac{12(\del-v)u}{u+1} \Big\} \na u\cdot\na v 
	- \io \frac{(u+1)^4}{(\del-v)^3} \cdot \Big\{ 2 - \frac{4(\del-v)u}{u+1}\Big\} |\na v|^2 \nn\\
	&\le& \io \frac{(u+1)^3}{(\del-v)^2} \cdot \big\{ 8 + 12\del \big\} |\na u| \, |\na v|
	- \frac{3}{2} \io \frac{(u+1)^4}{(\del-v)^3} |\na v|^2 \\
	&\le& \frac{1}{6} \io \frac{(u+1)^2}{\del-v} \cdot \big\{ 64 + 192\del  + 144\del^2 \big\} |\na u|^2 \\
	&=& 12 \io \frac{(u+1)^2}{\del-v} |\na u|^2
	- \io \frac{(u+1)^2}{\del-v} \cdot \Big\{ \frac{4}{3} - 32\del - 24\del^2 \Big\} |\na u|^2
	\qquad \mbox{for all } t>0.
  \eas
  Since here the number $c_1:=\frac{4}{3} - 32\del - 24\del^2$ is positive by (\ref{5.2}), from (\ref{5.3})
  we infer according to the nonpositivity of the rightmost summand therein that
  \bas
	\frac{d}{dt} \io \frac{(u-1)^4}{\del-v}
	+ c_1 \io \frac{(u+1)^2}{\del-v} |\na u|^2 
	\le 0
	\qquad \mbox{for all } t>0
  \eas
  and that thus, since trivially $\frac{(u+1)^2}{\del-v} \ge \del$,
  \bas
	\int_0^\infty \io |\na u|^2 < \infty.
  \eas
  According to mass conservation in (\ref{CC1}) and a Poincar\'e inequality, this implies that
  \bas
	\int_0^\infty \io |u-\ouz|^2 < \infty,
  \eas
  so that since straightforward application of parabolic theory turns the boundedness property in (\ref{bound}) 
  into uniform continuity of $u$ in $\bom\times [0,\infty)$, it readily follows that
  \be{5.4}
	u(\cdot,t)\to \ouz
	\quad \mbox{in } L^\infty(\Om)
	\qquad \mbox{as } t\to\infty.
  \ee
  As an integration in (\ref{CC1}) independently shows that
  \bas
	\int_0^\infty \io uv = \io v_0 < \infty,
  \eas
  using (\ref{5.4}) here we obtain that
  \bas
	\int_0^\infty \io v < \infty,
  \eas
  which completes the proof of (\ref{conv}), as also $v$ is uniformly continuous by parabolic H\"older regularity theory.
\qed
The rates of convergence in (\ref{conv}) have first been addressed in \cite{zhang_li} for a two-dimensional problem 
containing (\ref{CC1}) as a subsystem (cf.~also \cite{zhang_li_conv2}). 
The approach developed there, however, can readily be generalized so as to cover actually
all the cases from Theorems \ref{theo1} and \ref{theo3} and Proposition \ref{prop5}.
A simple but crucial observation in this regard uses the first part of the qualitative statement in (\ref{conv}) to make sure that
given $\eta>0$ one find $t_1(\eta)>0$ such that $u>\ouz-\eta$ and hence
\bas
	v_t \le \Del v - (\ouz-\eta) v
	\qquad \mbox{in } \Om\times (t_1(\eta),\infty).
\eas
By a comparison argument, this already implies an exponentially decaying upper bound for $\|v(\cdot,t)\|_{L^\infty(\Om)}$,
and suitably combining this with known smoothing properties of the Neumann heat semigroup and the mere boundedness of $u$ 
for large $t$ leads to an estimate of the form
\bas
	\|\na v(\cdot,t)\|_{L^\infty(\Om)} \le C(\eta) e^{-(\lam-\eta) t}
	\qquad \mbox{for all } t>t_2(\eta),
\eas
valid for each $\eta>0$ and some $C(\eta)>0$, with $\lam:=\min\{\ouz \, , \, \lam_1\}$ and $\lam_1>0$ denoting the smallest
nonzero eigenvalue of the Neumann Laplacian on $\Om$. 
This control of the taxis gradient facilitates a similar argument with respect to the first solution component, 
in summary leading to the following.
\begin{theo}\label{theo6} (\cite{zhang_li}) \quad
  Let $n\ge 1$ and $\Om\subset\R^n$ be a bounded domain with smooth boundary, and suppose that for some $T>0$,
  the pair $(u,v) \in \big( C^{2,1}(\bom\times [T,\infty))\big)^2$ solves the boundary value problem in (\ref{CC1})
  classically and is such that (\ref{conv}) holds. Then given any $\eta>0$ one can find $C(\eta)>0$ such that
  \be{6.1}
	\|u(\cdot,t)-\ouz\|_{L^\infty(\Om)} \le C(\eta) e^{-(\lam-\eta)t}
	\qquad \mbox{for all } t>T
  \ee
  and
  \be{6.2}
	\|v(\cdot,t)\|_{L^\infty(\Om)} \le C(\eta) e^{-(\ouz-\eta)t}
	\qquad \mbox{for all } t>T,
  \ee
  where $\ouz:=\frac{1}{|\Om|} \io u(\cdot,T)$, and where $\lam:=\min\{\ouz \, , \, \lam_1\}$, with $\lam_1>0$ denoting
  the principal eigenvalue of $-\Del$ under homogeneous Neumann boundary conditions on $\Om$.
\end{theo}
Extensions of this to certain strip-like domains are discussed in \cite{wu_lin}.\abs 
In some application contexts, chemotactic motion cannot be assumed to be adequately described by cross-diffusion operators
as simple as that in (\ref{CC1}). Typical examples for considerable deviations include bacterial chemotaxis near surfaces
which, in fact, may contain rotational components orthogonal to the signal gradient (\cite{xue_othmer}), or also
non-negligible dependence of chemotactic responses on the concentration of the respective chemical (\cite{goldstein2005}). 
As a generalization of (\ref{CC1}) accounting for both these possibilities, let us consider the problem
\be{CC2}
	\lball
	u_t=\Del u - \na\cdot \big( uS(x,u,v)\cdot\na v\big),
	\qquad & x\in\Om, \ t>0, \\[1mm]
	v_t = \Del v - uv
	\qquad & x\in\Om, \ t>0, \\[1mm]
	\big( uS(x,u,v)\cdot\na v\big) \cdot \nu = \na v\cdot \nu=0,		
	\qquad & x\in\pO, \ t>0, \\[1mm]
	u(x,0)=u_0(x), \quad v(x,0)=v_0(x),
	\qquad & x\in\Om,
	\ear
	\tag{CC2}
\ee
where now $S$ is a suitably regular function that attains values in $\R^{n\times n}$;
as prototypical choices, one may think of tensor-valued $S$ containing off-diagonal entries such as, in the case $n=2$,
\bas
	S(x,s,\sig):= 
	\chi \cdot \left( \begin{array}{cc}
	1 & 0 \\ 0 & 1
	\end{array} \right)
	+ \beta \cdot \left( \begin{array}{cc}
	0 & 1 \\ -1 & 0
	\end{array} \right),
	\qquad (x,s,\sig)\in \bom\times [0,\infty)^2,
\eas
or also
\bas
	S(x,s,\sig):=\Theta(\sig-\sig_\star) {\bf E}_n
	\qquad (x,s,\sig)\in \bom\times [0,\infty)^2
\eas
with some $\sig_\star>0$, where ${\bf E}_n$ denotes the $n\times n$ unit matrix, and where $\Theta$ is some suitably smooth approximation of the Heaviside function (\cite{xue_othmer}, \cite{goldstein2005}).\abs
In such more complex settings, subtle structures as that expressed in (\ref{en}) apparently can no longer be expected. 
Consequently, the question how far taxis-driven blow-up is precluded by dissipation again becomes nontrivial
already in two-dimensional domains in which (\ref{CC1}) allowed for the comprehensive result from Theorem \ref{theo1}.
In the following we briefly reproduce an approach by which it becomes possible to nevertheless provide at least some basic
indication for the guess that the consumption mechanism in (\ref{CC2}) should exert some relaxing
influence in comparison to situations in Keller-Segel-production systems.\abs
In the absence of any further assumptions on $S$ other than regularity, it seems that the only expedient 
option to handle taxis is to simply estimate its effect in modulus.  
Here a favorable difference between (\ref{CC2}) and (\ref{KS}) consists in the fact that solutions of arbitrary size,
and in domains of any dimension, should satisfy
\be{gradv}
	\int_0^T \io |\na v|^2 \le \frac{1}{2} \io v_0^2
	\qquad \mbox{for all } T>0.
\ee
Though rather poor with respect to its topological framework, this unconditional bound for the cross-diffusive gradient
can indeed be used to derive some basic regularity feature of the key quantity $u$: Namely, analyzing the evolution
properties of $\io \ln (u+1)$ readily shows that when $S$ is merely assumed to be bounded, estimates of the form
\be{gradlogu}
	\int_0^T \io \frac{|\na u|^2}{(u+1)^2} \le C
\ee
should be available for all $T>0$.\abs
Following this strategy in the course of an approximation by global smooth solutions to regularized variants of (\ref{CC2}),
one can indeed derive a result on global existence, albeit
within a solution concept yet considerably weaker than that from Theorem \ref{theo2}.
\begin{theo}\label{theo7} (\cite{win_SIMA2015}) \quad
  Let $n\ge 2$ and $\Omega\subset\R^n$ be a bounded domain with smooth boundary, let 
  $S\in C^2(\bom\times [0,\infty)\times [0,\infty);ℝ^{n\times n})$ be bounded, and assume (\ref{init}).   
  Then there exists at least one pair $(u,v)$ of nonnegative functions fulfilling (\ref{reg2}) as well as
  \be{reg_log}
	\ln (u+1) \in L^2_{loc}([0,\infty);W^{1,2}(\Omega)), 
  \ee
  such that (\ref{CC2}) is solved in the sense that (\ref{2.3}) is satisfied for any $\varphi\in C_0^\infty(\bom\times [0,\infty))$,
  that
  \be{7.1}
	\io u(\cdot,t) \le \io u_0
	\qquad \mbox{for a.e.~} t>0,
  \ee
  and that for each nonnegative $\varphi \in C_0^\infty(\bar\Omega \times [0,\infty))$,
  the inequality 
  \bea{7.2}
	- \int_0^\infty \io \ln (u+1) \varphi_t 
	&-&\io \ln(u_0+1) \varphi(\cdot,0)\nn\\
	&\ge& \int_0^\infty \io |\nabla \ln (u+1)|^2 \varphi 
	- \int_0^\infty \io \na \ln(u+1) \cdot\na\varphi \nn\\
	& & 
	- \int_0^\infty \io \frac{u}{u+1} \nabla \ln (u+1) \cdot \Big(S(x,u,v) \cdot \nabla v \Big) \varphi \nn\\
	& &
	+ \int_0^\infty \io \frac{u}{u+1} \Big(S(x,u,v) \cdot\nabla v \Big) \cdot\nabla \varphi 
  \eea
  holds.
\end{theo}
The problem of determining how far solutions to (\ref{CC2}) may enjoy regularity properties beyond those from Theorem \ref{theo7}
forms a widely open challenge. Only in planar cases,
some 
results on eventual smoothness and large time
stabilization in the style of Theorem \ref{theo3} are available under the general assumptions on $S$ from Theorem \ref{theo7},
provided that $\Om$ is convex (\cite{heihoff}, \cite{win_IMRN}).
In higher-dimensional situations, it seems in general even unclear whether the natural mass conservation property $\io u=\io u_0$
can be expected instead of the inequality in (\ref{7.1}).
After all, both this and also eventual regularity as well as (\ref{conv}) could at least 
be asserted in frameworks of radially symmetric solutions to (\ref{CC2}) in balls of arbitrary dimension, assuming that
$S=S(u,v)$ is sufficiently smooth and bounded (\cite{weirun_tao}).\abs
An important observation now confirms that not only the result on global existence of classical small-signal solutions
from Theorem \ref{theo4} extends to general bounded matrix-valued $S$, but that moreover also any algebraic-type saturation
of cross-diffusive fluxes at large population densities prevents blow-up. 
This is achieved in \cite{zhang_MANA2016} by analyzing a functional of the form $\F_2$ from (\ref{F2}) with
a choice of $\vp$ similar to that in (\ref{phi}) (see also \cite{cao_wang_DCDSB2015} for a precendent):
\begin{theo}\label{theo8} (\cite{zhang_MANA2016}) \quad
  Let $n\ge 2$ and $\Omega\subset\R^n$ be a bounded domain with smooth boundary, assume (\ref{init}), and let 
  $S\in C^2(\bom\times [0,\infty)\times [0,\infty);ℝ^{n\times n})$ be bounded
  and such that either
  \be{S0}
	|S(x,s,\sig)| \le S_0(\sig)
	\quad \mbox{for all } (x,s,\sig)\in \bom\times [0,\infty)^2
	\qquad \mbox{with some nondecreasing $S_0:[0,\infty)\to [0,\infty)$}
  \ee
  fulfilling
  \bas
	S_0\big( \|v_0\|_{L^\infty(\Om)} \big) < \frac{2n}{3n(11n+2)},
  \eas
  or, alternatively, 
  \bas
	|S(x,s,\sig)| \le \frac{C}{(s+1)^\theta}
	\qquad \mbox{for all $x\in\bom, usge 0$ and $\sig\in [0,\|v_0\|_{L^\infty(\Om)}]$ and some } C>0.
  \eas
  Then (\ref{CC2}) admits a global bounded classical solution satisfying (\ref{reg}) as well as (\ref{conv}).
\end{theo}
\mysection{Including nonlinear diffusion and singular sensitivities}
Various types of nonlinear diffusion mechanisms play important roles in microbial motion (\cite{kowalczyk},
\cite{hillen_painter_VFE}, \cite{hillen_painter_survey}),
and corresponding generalizations of (\ref{CC1}) and (\ref{CC2}), such as
\be{CC3}
	\lball
	u_t=\na \cdot \big( D(u) \na u\big) u - \na\cdot \big(uS(x,u,v)\cdot \na v\big), 
	\qquad & x\in\Om, \ t>0, \\[1mm]
	v_t = \Del v - uv
	\qquad & x\in\Om, \ t>0, \\[1mm]
	\big( D(u) \na u\big) u - uS(x,u,v)\cdot \na v \big)\cdot \nu = \na v\cdot\nu=0,
	\qquad & x\in\pO, \ t>0, \\[1mm]
	u(x,0)=u_0(x), \quad v(x,0)=v_0(x),
	\qquad & x\in\Om,
	\ear
	\tag{CC3}
\ee
have quite intensely been studied in the literature. 
Especially in light of the standoff with regard to global smooth solvability of (\ref{CC1}) in three-dimensional domains, 
but also of (\ref{CC2}) in planar situations for general bounded $S$,
a particular focus in this context is how far relaxing effects can be expected to result from a presupposed
enhancement of diffusion at large population densities, such as present when
in the prototypical choice of porous medium type diffusion, that is, in
\be{PME}
	D(s)=s^{m-1},
	\qquad s\ge 0,
\ee
the adiabatic exponent is assumed to satisfy $m>1$.
Due to the degeneracy at vanishing population densities which additionally goes along with such a simple
functional setting, classical solutions
can apparently not be expected to exist for general nonnegative initial data. 
Accordingly, a focus is commonly set on the construction of global weak solutions with uniformly bounded components.\abs
Indeed, in two-dimensional convex domains and for arbitrary bounded smooth matrix-valued functions $S$, 
such bounded weak solutions exist under the mere assumption that $m>1$.
The following result in this regard extends previous partial findings obtained in \cite{DiFLM}.	
\begin{theo}\label{theo9} (\cite{cao_ishida}) \quad
  Let $\Omega\subset\R^2$ be a bounded convex domain with smooth boundary, assume (\ref{init}), let 
  $S\in C^2(\bom\times [0,\infty)\times [0,\infty);ℝ^{n\times n})$ be such that (\ref{S0}) holds, and suppose that
  (\ref{PME}) is valid
  with some
  \bas
	m>1.
  \eas
  Then there exist nonnegative functions
  \be{reg3}
	\lbal
	u\in L^\infty(\Om\times (0,\infty))
	\qquad \mbox{and} \\[1mm]
	v\in L^\infty((0,\infty);W^{1,\infty}(\Om))
	\ear
  \ee
  such that 
  \be{9.1}
	\na u^m \in L^2_{loc}(\bom\times [0,\infty);\R^2)
  \ee
  and that (\ref{CC3}) is solved in the sense that
  \be{9.2}
	- \int_0^\infty \io u\vp_t - \io u_0\vp(\cdot,0)
	= - \int_0^\infty \io \na u^m\cdot\na \vp
	+ \int_0^\infty \io u \big(S(x,u,v)\cdot \na v\big) \cdot\na\vp
  \ee
  as well as (\ref{2.3})
  hold for all $\vp\in C_0^\infty(\bom\times [0,\infty))$.
\end{theo}
The derivation of a three-dimensional analogue has apparently been forming a considerable challenge for almost a decade, with
partial results having successively been improved in several steps;
only recently, two independent approaches have led to an essentially complete picture in this regard.
One of these uses an adaptation of the energy inequality (\ref{en}) to the present setting
involving nonlinear diffusion, and uses corresponding basic estimates in order to
identify certain quasi-energy properties of
\bas
	\frac{d}{dt} \io \Big\{ \frac{|\na u|^2}{u} + |\Del v|^2 + |\na v|^4 \Big\}.
\eas
In \cite{jin},
this indeed leads to a result parallel to that from Theorem \ref{theo9} under stronger assumptions on the initial data,
and under the additional condition that $S={\bf E}_n$ is scalar.\abs
By means of second and independent approach, it becomes possible not only to remove these restrictions, 
but also to actually cover diffusivities with arbitrarily slow growth.
In fact, by making essential use of some favorable evolution properties of the mixed functional 
$\io \frac{u}{v} |\na v|^2$ it is possible to see that whenever $D$ is smooth on $[0,\infty)$ and positive on $(0,\infty)$
with $D(s)\to + \infty$ as $s\to\infty$, and whenever the matrix-valued $S$ is bounded,
for some suitably chosen smooth nonnegative $\psi$ on $[0,\infty)$ an inequality of the form
\bas
	\frac{d}{dt} \io \Big\{ \psi(u) + d_1 \frac{u}{v} |\na v|^2 + d_2 \frac{1}{v^3} |\na v|^4 \Big\}
	+ a \io \Big\{ \psi(u) + d_1 \frac{u}{v} |\na v|^2 + d_2 \frac{1}{v^3} |\na v|^4 \Big\}
	+ \io |\na u|^2
	\le b
\eas
holds for $t>0$ with some positive $d_1, d_2, a$ and $b$.\abs
Conclusions drawn from this can be used, even under slightly milder assumptions on $S$ than mentioned above, 
to complete a considerable history
on blow-up prevention by arbitrarily mild 
diffusion enhancement in the three-dimensional version of (\ref{CC3}),
as actually developed by several authors to be mentioned here:
\begin{theo}\label{theo10} (\cite{duan_xiang}, \cite{kang_JMP}, \cite{chung_kang_kim2014}, \cite{jin}, 
\cite[Corollary 1.2 and p.95]{win_DCDSB}) \quad
  Let $\Om\subset\R^3$ be a bounded domain with smooth boundary, and suppose that
  \be{S}
	\lbal
	S\in C^2(\bom\times [0,\infty)\times (0,\infty);\R^{3\times 3})
	\mbox{ is such that } \\[2mm]
	\displaystyle |S(x,s,\sig)| \le \frac{S_0(\sig)}{\sig^\frac{1}{2}} \qquad \mbox{ for all $(x,s,\sig)\in\Om\times (0,\infty)^2$}  \\[2mm]
	\hs{41mm} \mbox{with some nondecreasing $S_0: (0,\infty)\to (0,\infty)$,}
	\ear
  \ee
  and that either (\ref{PME}) holds with some $m>1$, or that
  \bas
	D\in \bigcup_{\vartheta\in (0,1)} C^\vartheta_{loc}([0,\infty)) \cap C^2((0,\infty))
	\qquad \mbox{is positive on } (0,\infty)
  \eas
  and such that $\liminf_{s\searrow 0} \frac{D(s)}{s}>0$ as well as
  \be{10.1}
	D(s) \to +\infty
	\qquad \mbox{as } s\to\infty.
  \ee
  Then given any initial data fulfilling (\ref{init}),
  one can find functions $u$ and $v$ which are such that (\ref{reg3}) holds, that $u\ge 0$ and $v>0$ a.e.~in $\Om\times (0,\infty)$
  and $u|S(\cdot,u,v)|\in L^1_{loc}(\bom\times [0,\infty))$, 
  that (\ref{2.3}) holds for all $\vp\in C_0^\infty(\bom\times [0,\infty))$,
  and that writing $D_1(s):=\int_0^s D(\sig) d\sig$ for $s\ge 0$
  we have
  \be{10.2}
	- \int_0^\infty \io u\vp_t - \io u_0\vp(\cdot,0)
	= \int_0^\infty \io D_1(u) \Del\vp
	+ \int_0^\infty \io u \big(S(x,u,v)\cdot \na v\big) \cdot\na\vp
  \ee
  for all $\vp\in C_0^\infty(\bom\times [0,\infty))$ satisfying $\frac{\pa\vp}{\pa\nu}=0$ on $\pO\times (0,\infty)$.
%
\end{theo}
In domains of dimension $n\ge 4$, more restrictive assumptions on $D$ seem necessary to ensure global existence of
bounded solutions, see \cite{yilong_wang_xiang}, \cite{fan_jin}, \cite{liangchen_wang2014} and \cite{liangchen_wang2015}, for instance;
the to date mildest assumption on $m$ in (\ref{CC3})-(\ref{PME}) for such a conclusion to hold for $S={\bf E}_n$
seems to require that $m>\frac{3n}{2n+2}$ (\cite{jiashan_zheng_yifu_wang2017}). Qualitative properties such as stabilization toward semitrivial steady states and propagation features of $\supp u$ have been discussed in \cite{jin_wang_yin} and \cite{xu_ji_mei_yin}.\abs
Combined effects of nonlinear porous medium type diffusion and sensitivities $S=S(x,u,v)$ either growing or decreasing 
in an essentially algebraic manner as $u\to\infty$ have been studied in \cite{jiashan_zheng2022}, \cite{jiashan_zheng_JDE2022},
\cite{yilong_wang1}, \cite{yilong_wang2}, \cite{wang_mu_hu}, \cite{wang_zhang_zhang}
and \cite{wang_xie}, for instance.
Interesting examples addressing exponentially decaying diffusion rates can be found in \cite{liangchen_wang_CMA2017}
and \cite{bingchen_liu}.\\[5mm]
To appropriately cope with situations in which the intensity of chemotactic response depends on the signal in a monotonically 
decreasing manner, Keller and Segel \cite{keller_segel_traveling} proposed the system
\be{CC4}
	\lball
	u_t=\Del u - \na\cdot \big(uS(v)\na v\big), 
	\qquad & x\in\Om, \ t>0, \\[1mm]
	v_t = \Del v - uv
	\qquad & x\in\Om, \ t>0, \\[1mm]
	\big( \na u - uS(v)\na v \big)\cdot \nu = \na v\cdot\nu=0,
	\qquad & x\in\pO, \ t>0, \\[1mm]
	u(x,0)=u_0(x), \quad v(x,0)=v_0(x),
	\qquad & x\in\Om,
	\ear
	\tag{CC4}
\ee
with the particular choice
\be{sing}
	S(s)=\frac{1}{s}, 
	\qquad s>0,
\ee
corresponding to the Weber-Fechner law of stimulus response (\cite{rosen}).
In the presence of such singular sensitivities, this specific version of (\ref{CC3}) appears to be mathematically quite
delicate, as witnessed by the circumstance that although several findings indicate a strong support of wave-like
structured behavior (see \cite{kang_waves} and also \cite{zhian_wang_DCDS} for a survey),
already in two-dimensional cases neither a comprehensive theory
of global classical solvability nor any result on the occurrence of blow-up phenomena seem to have been established yet.\abs
After all, globally smooth solutions can be constructed at least for suitably small initial data. Indeed, in \cite{black_JDE2018}
it was observed that for arbitrary $\mu>0$, the expression
\bas
	\F_3(t):= \io u\ln \frac{u}{\mu} + \frac{1}{2} \io \frac{|\na v|^2}{v^2}
\eas
acts as a conditional energy functional for the planar version of (\ref{CC4}) in the sense that smooth solutions satisfy
\be{en3}
	\F_3'(t) + \io \frac{|\na u|^2}{u}
	+ \bigg\{ \frac{1}{2} - a \io \frac{|\na v|^2}{v^2} - a \io u_0 \bigg\} \cdot \io |\Del \ln v|^2 \le 0
\ee
with some $a=a(\Om)$. Suitably relating the ill-signed part on the right to $\F_3$ itself leads to 
bounds for $\F_3$ and the dissipation rate functionals $\io \frac{|\na u|^2}{u}$ as well as $\io |\Del\ln v|^2$
which in planar cases can be used as a starting point for a bootstrap-type series of regularity arguments,
finally leading to the following result on global existence of small-data solutions:
\begin{theo}\label{theo11} (\cite{black_JDE2018}) \quad
  Suppose that $\Om\subset\R^2$ is a bounded domain with smooth boundary, and assume (\ref{sing}).
  Then one can find $K>0$ and $m>0$ with the property that whenever $u_0$ and $v_0$ satisfy (\ref{init}) and are such that
  \bas
	\io u_0<m\quad\text{and}\quad\io u_0 \ln \frac{u_0}{\mu} + \frac{1}{2} \io \frac{|\na v_0|^2}{v_0^2} \le K - \frac{\mu |\Om|}{e}
  \eas
  for some $\mu>0$, 
  there exist functions $u$ and $v$ which are such that (\ref{reg}) holds, that $u\ge 0$ and $v>0$ in $\bom\times [0,\infty)$,
  and that (\ref{CC4}) is solved classically.
  Moreover, this solution stabilizes in the sense that (\ref{conv}) is valid.
\end{theo}
For large initial data, only some generalized solutions have been found to exist. 
The following result in this direction can be established by means of a far relative of the reasoning near Theorem \ref{theo7},
supplemented by an application of the two-dimensional Moser-Trudinger inequality in deriving certain $L^1$ compactness
properties in the first solution component.
Accordingly, our statement here goes beyond that from Theorem \ref{theo7} in asserting that the obtained solutions
indeed enjoy the natural mass conservation property formally associated with (\ref{CC4}).
\begin{theo}\label{theo12} (\cite{wang16}, \cite{black_JDE2018}, \cite{win_M3AS}) \quad
  Let $\Omega\subset\R^2$ be a bounded domain with smooth boundary, and assume (\ref{sing}).
  Then for all $u_0$ and $v_0$ satisfying (\ref{init}), the problem (\ref{CC4}) admits at least one 
  global generalized solution in the sense that there exist functions $u$ and $v$ such that (\ref{reg2}) holds, that
  $u\ge 0$ and $v>0$ a.e.~in $\Om\times (0,\infty)$, 
  \bas
	\ln (u+1) \in L^2_{loc}([0,\infty);W^{1,2}(\Omega))
	\qquad \mbox{and} \qquad
	\ln v \in L^2_{loc}([0,\infty);W^{1,2}(\Omega)),
  \eas
  that
  \bas
	\io u(\cdot,t) = \io u_0
	\qquad \mbox{for a.e.~} t>0,
  \eas
  and that for each nonnegative $\varphi \in C_0^\infty(\bar\Omega \times [0,\infty))$, the identity in (\ref{2.3}) as well as
  the inequality
  \bas
	- \int_0^\infty \io \ln (u+1) \varphi_t 
	&-&\io \ln(u_0+1) \varphi(\cdot,0) \\
	&\ge& \int_0^\infty \io |\nabla \ln (u+1)|^2 \varphi 
	- \int_0^\infty \io \na \ln(u+1) \cdot \na \varphi \\
	& & 
	- \int_0^\infty \io \frac{u}{u+1} \big(\nabla \ln (u+1) \cdot \na \ln v \big) \varphi \\
	& &
	+ \int_0^\infty \io \frac{u}{u+1} \na \ln v \cdot \nabla \varphi 
  \eas
  hold.
\end{theo}
Now the observation that solutions to suitably regularized variants of (\ref{CC4}) satisfy
\bas
	\frac{1}{T} \int_0^T \io \frac{|\na v|^2}{v^2} \le \frac{1}{T} \io \ln \frac{\|v_0\|_{L^\infty(\Om)}}{v_0} + \io u_0
	\qquad \mbox{for all } T>0,
\eas
facilitates a fruitful application of (\ref{en3}) also to trajectories more general than those from Theorem \ref{theo11}.
In fact, under a smallness assumption merely involving a total mass functional as a quantity of immediate biological relevance,
the following conclusion on eventual regularity can thereby be drawn.
\begin{theo}\label{theo13} (\cite{black_JDE2018}, \cite{ji_liu_JEE2021}) \quad
  Let $\Om\subset\R^2$ be a bounded domain with smooth boundary, and assume (\ref{sing}).
  Then there exists $M>0$ with the property that whenever (\ref{init}) holds with
  \bas
	\io u_0 \le M,
  \eas
  one can find $T>0$ such that the global generalized solution $(u,v)$ of (\ref{CC4}) from Theorem \ref{theo13} satisfies
  $(u,v)\in (C^{2,1}(\bom\times [T,\infty)))^2$ and $v>0$ in $\bom\times [T,\infty)$. Moreover, (\ref{conv}) holds, and additionally
  we have
  \bas
	\frac{\na v(\cdot,t)}{v(\cdot,t)} \to 0
	\quad \mbox{in } L^\infty(\Om)
	\qquad \mbox{as } t\to\infty.
  \eas
\end{theo}
In higher-dimensional versions of (\ref{CC4}), only in radially symmetric settings some global solutions seem to have been 
constructed so far, and only in strongly generalized frameworks of renormalized solutions (\cite{win_renormalized}).
For a Cauchy problem in $\R^n$ associated with (\ref{CC4}), under some conditions inter alia
requiring smallness of both $\|u_0\|_{W^{1,2}(\R^n)}$ and $\|\ln v_0\|_{W^{1,2}(\R^n)}$, a result on global existence 
and large time stabilization of certain strong solutions was derived in \cite{wang_xiang_yu}.\abs
Signal-dependent sensitivities less singular than those in (\ref{sing}) have been studied in contexts of essentially
algebraic behavior, as prototypically modeled by choices of the form
\be{sing1}
	S(s)=\frac{1}{s^\al},
	\qquad s>0,
\ee
with 		
$\al\in (0,1)$. Two results in this direction are summarized in the following.
\begin{theo}\label{theo14} (\cite{kang_MMAS2021}, \cite{viglialoro}) \quad
  Let $\Om\subset\R^2$ be a bounded domain with smooth boundary, and assume (\ref{sing1}) with $\al>0$.\abs
%
  i) \ If 
  \bas
	\alpha<-\frac{1}{2}+\frac{1}{2}\cdot \Big(\frac{8\sqrt{2}+11}{7}\Big)^\frac{1}{2}\approx 0.3927,
  \eas
  then for arbitrary $u_0$ and $v_0$ complying with (\ref{init}), the problem (\ref{CC4}) admits a global 
  classical solution satisfying (\ref{reg}) as well as $v>0$ in $\bom\times [0,\infty)$.\abs
  ii) \ If $\al<1$, then there exist $δ>0$ and $m>0$ such that whenever (\ref{init}) is valid with 
  \bas
	\|v_0\|_{L^\infty(\Om)} \le \del \quad \text{and} \quad \io u_0<m,
  \eas
  one can find a global classical solution of (\ref{CC4}) for which (\ref{reg}) holds and $v>0$
  in $\bom\times [0,\infty)$.
\end{theo}
Beyond this, a statement on boundedness of the solutions from Theorem \ref{theo14} ii) is available 
under the assumption that, with some $\eta=\eta(\al)>0$,
\bas
	\|u_0\|_{L^1(\Om)}^\frac{2}{n} \cdot \Big\{ \|v_0\|_{L^\infty(\Om)} + \|v_0\|_{L^\infty(\Om)}^\theta \Big\} \le \eta
\eas
holds (\cite{win_DCDSB}).\abs
If the first equation in (\ref{CC4}) is replaced by
\bas
	u_t=\Del u - \na\cdot \Big( \frac{f(u)}{v} \na v\Big),
\eas
with the nonnegative function $f$ reflecting suitably strong saturation effects in the sense that 
\bas
	f(s) \le C (s+1)^{1-\frac{n}{4}-\eta}
	\quad \mbox{for all } s\ge 1
	\qquad \mbox{with some $C>0$ and $\eta>0$,}
\eas
global classical solutions can be found in $n$-dimensional domains with $n\ge 2$ (\cite{dongmei_liu}).
If, instead, (\ref{CC4}) is considered with its second equation altered to
\bas
	v_t=\Del v - g(u) v,
\eas
where $0\le g\in C^1([0,\infty)$ satisfies $g(0)=0$ and 
\bas
	g(s) \le C (s+1)^{1-\eta}
	\quad \mbox{for all } s\ge 1
	\qquad \mbox{with some $C>0$ and $\eta>0$,}
\eas
then in the corresponding two-dimensional boundary value problem global classical solutions exist (\cite{lankeit_viglialoro}).\abs
Some further indications for the guess that
strong singularities such as those in (\ref{sing}) may go along with a significant loss of regularity can be gained
upon observing that in contrast to the situation in (\ref{CC2}), for corresponding quasilinear variants of (\ref{CC4})
the literature so far provides results on global
existence of solutions which are at least locally bounded only under considerably more restrictive assumptions on enhancement
of diffusion. Specifically, if the first equation in (\ref{CC4}) is changed so as to become
\bas
	u_t=\na \cdot \big(D(u)\na u\big) - \na \cdot \Big(\frac{u}{v}\na v\Big)
\eas
with $D$ being suitably smooth and satisfying $D(s)\ge C s^{m-1}$ for all $s\ge 0$ with some $C>0$,
then for such locally bounded global weak solutions to exist in $n$-dimansional domains, 
to date the literature requires the assumption that $m>1+\frac{n}{4}$ (\cite{lankeit_singularconsumption})
which already in the case $n=2$ 
is considerably stronger than the corresponding condition $m>1$ underlying the boundedness result in Theorem \ref{theo9}.\abs
Findings concerned with further related variants of (\ref{CC4}) are contained in \cite{zhao}, \cite{viglialoro_AMO},
\cite{jia_yang} and \cite{jia_yang2}, \cite{yan_li_yuxiang_sing_cons}, for instance.
\mysection{Couplings to further mechanisms}		
The intention of this section is to review some developments related to
the embedding of chemotaxis-consumption systems into more complex models involving further interaction mechanisms and components.\abs
A first example in this regard is provided by the class of logistic Keller-Segel-absorption systems,
\be{CC5}
	\lball
	u_t=\Del u - \na\cdot \big(uS(x,u,v)\cdot \na v\big) + f(u), 
	\qquad & x\in\Om, \ t>0, \\[1mm]
	v_t = \Del v - uv
	\qquad & x\in\Om, \ t>0, \\[1mm]
	\frac{\pa u}{\pa\nu}=\frac{\pa v}{\pa\nu}=0,
	\qquad & x\in\pO, \ t>0, \\[1mm]
	u(x,0)=u_0(x), \quad v(x,0)=v_0(x),
	\qquad & x\in\Om,
	\ear
	\tag{CC5}
\ee
again with matrix-valued $S$, and 
with $f$ representing zero-order mechanisms of proliferation and death; 
typical choices include
\be{log}
	f(s) = \rho s - \mu s^\kappa,
	\qquad s\ge 0,
\ee
with $\rho\in\R$, $\mu>0$ and $\kappa>1$.\abs
Now a first effect of the superlinear damping thereby introduced becomes manifest in the fact that beyond
global generalized solvability, also in the presence of widely arbitrary matrices $S$ 
mass conservation can be asserted, as implicitly contained in the proof
presented in \cite{wenbin_JDE2021} (cf.~also \cite{wenbin_NLA}, where tensor-valued sensitivities are explicitly
included).
\begin{theo}\label{theo15} (\cite{wenbin_JDE2021}) \quad
  Let $n\ge 2$ and $\Omega\subset\R^n$ be a bounded domain with smooth boundary, let 
  $S\in C^2(\bom\times [0,\infty)\times [0,\infty);ℝ^{n\times n})$ be bounded, and let $f$ be as in (\ref{log})
  with $\rho\in\R$, $\mu>0$ and $\kappa>1$.
  Then for any choice of $u_0$ and $v_0$ satisfying (\ref{init}), it is possible to find nonnegative functions
  $u$ and $v$ which are such that (\ref{reg2}) and (\ref{reg_log}) holds, that
  \bas
	u \in L^\kappa_{loc}(\bom\times [0,\infty))
  \eas
  with
  \bas
	\io u(\cdot,t) = \io u_0
	\qquad \mbox{for a.e.~} t>0,
  \eas
  and that
  whenever $\varphi\in C_0^\infty(\bom\times [0,\infty))$, (\ref{2.3}) as well as the inequality
  \bas
	- \int_0^\infty \io \ln (u+1) \varphi_t 
	&-&\io \ln(u_0+1) \varphi(\cdot,0)\nn\\
	&\ge& \int_0^\infty \io |\nabla \ln (u+1)|^2 \varphi 
	- \int_0^\infty \io \na \ln(u+1) \cdot\na\varphi \\
	& & 
	- \int_0^\infty \io \frac{u}{u+1} \nabla \ln (u+1) \cdot \Big(S(x,u,v) \cdot \nabla v \Big) \varphi \nn\\
	& &
	+ \int_0^\infty \io \frac{u}{u+1} \Big(S(x,u,v) \cdot\nabla v \Big) \cdot\nabla \varphi 
	+ \int_0^\infty \io \frac{f(u)}{u+1} \vp
  \eas
  are valid.
\end{theo}
Interestingly, the most prototypical choice $\kappa=2$, as most commonly used in applications,
seems to play also a technical role, since it marks the minimal logistic degradation effect for which global
classical solutions can be found. In fact, the following statement in this regard concentrates
on scalar $S$, but could readily be extended to bounded smooth matrix-valued $S$ as well:
\begin{theo}\label{theo16} (\cite{lankeit_wang}) \quad
  Let $n\ge 2$ and $\Omega\subset\R^n$ be a bounded domain with smooth boundary,
  and let $f$ be as in (\ref{log}) with $\rho\in\R$, $\mu>0$ and 
  \bas
	\kappa=2.
  \eas
  Then there exists $C>0$ such that whenever (\ref{init}) holds with
  \bas
	\mu> C\|v_0\|_{L^\infty(\Om)}^\frac{1}{n} + C\|v_0\|_{L^\infty(\Om)}^{2n},
  \eas
  the problem (\ref{CC5}) admits a global classical solution satisfying (\ref{reg}) as well as
  $\sup_{t>0} \|u(\cdot,t)\|_{L^\infty(\Om)} <\infty$.\abs
  If, moreover, $\rho>0$, then
  \be{conv1}
	u(\cdot,t) \to \frac{\rho}{\mu}
	\quad \mbox{and} \quad
	v(\cdot,t) \to 0
	\quad \mbox{in } L^\infty(\Om)
	\qquad \mbox{as } t\to\infty.
  \ee
\end{theo}
For the generalized solutions from Theorem \ref{theo15}, similar statements on asymptotic stabilization are available
actually for arbitrary $\kappa>1$,
albeit partially in weaker topological settings (\cite{wenbin_NLA}, \cite{cao_zhuang});
in three-dimensional versions of (\ref{CC5}) with scalar $S$, these solutions are even known to become smooth
eventually (\cite{cao_zhuang}; cf.~also \cite{yulan_M3AS}).\abs
For (\ref{CC5}) and some quasilinear variants 
involving singular sensitivities, some results on existence and qualitative properties have been documented in 
\cite{wei_wang}, 
\cite{lankeit_lankeit1}, \cite{lankeit_lankeit2}, \cite{zhao_zheng_logistic}, \cite{pang_yifu_wang}, \cite{he}, \cite{wenji},
\cite{jiashan_zheng11} and \cite{chunhua_jin2}.\abs
A second example of more involved mechanisms for interaction is given by the following equations for indirect signal absorption where the signal substance is not directly consumed upon contact with the motile cells but instead first an intermediate component (of concentration $w$) is produced which then reduces the signal: 
\be{CC55}
\lball 
 u_t = \Del u - \na \cdot(u\na v)\qquad & x\in\Om, \ t>0, \\[1mm]
 v_t = \Del v - vw \qquad & x\in \Om, \ t>0,\\[1mm]
 w_t = -\delta w + u\qquad & x\in \Om, \ t>0,\\[1mm]
 	\frac{\pa u}{\pa\nu}=\frac{\pa v}{\pa\nu}=0,
	\qquad & x\in\pO, \ t>0, \\[1mm]
	u(x,0)=u_0(x), \quad v(x,0)=v_0(x), \quad w(x,0)=w_0(x)
	\qquad & x\in\Om,
 	\ear
 	\tag{CC6}
\ee
While having exchanged $u$ for $w$ in the second equation voids the cancellation on which \eqref{en} relied, a reasoning based on \eqref{F2} again makes it possible to obtain global small-data solutions (with the condition from \cite{fuest_indirectabsorption} improved in \cite{frassu_viglialoro} and \cite{liu_li_huang}, similarly as in Theorem~\ref{theo4}).  

Supported by the observation that spatio-temporal $L^p$-bounds for $u$ entail temporally uniform $L^p(\Omega)$ bounds for $w$, a study of $\mathcal{F}_2$ from \eqref{F2}, now for $p\in(0,1)$ and a new suitable choice of $\varphi$, reveals global existence of classical solutions (\cite{fuest_indirectabsorption}).  

Extensions involving logistic source terms (\cite{liu_li_huang}) or nonlinear diffusion and sensitivity (\cite{zhang_liu,zheng_xing_20}) are available, as are studies of \eqref{CC55} with additional diffusion of the intermediate (e.g. \cite{zheng_xing_22,xing_zheng_xiang_wang_21}). 

Another class of systems which can be viewed as extensions of (\ref{CC1}) is formed by the prey-taxis models (\cite{bendahmane})
\be{CC6}
	\lball
	u_t = \Del u - \chi \na \cdot (u\na v)		
		+ \gamma uh(v) + f(u)
	\qquad & x\in\Om, \ t>0, \\[1mm]
	v_t = \Del v - uh(v) + g(v)
	\qquad & x\in\Om, \ t>0, \\[1mm]
	\frac{\pa u}{\pa\nu}=\frac{\pa v}{\pa\nu}=0,
	\qquad & x\in\pO, \ t>0, \\[1mm]
	u(x,0)=u_0(x), \quad v(x,0)=v_0(x),
	\qquad & x\in\Om,
	\ear
	\tag{CC7}
\ee
and their close relatives.
Here, the evolution of $v$ is now influenced not only by the absorption mechanism, expressed here through the
summand $-uh(v)$ in a form slightly more gerenal than in the above systems, but moreover potentially also by the term $+g(v)$
representing reproduction in the considered prey population density $v$.
The corresponding population density $u$ of predators, in turn, benefits from consumption of preys through the term $+\gamma uh(v)$,
and additionally may be affected by spontaneous proliferation and death according to a law for $f$ preferably of the style
in (\ref{log}).\abs
Under appropriate technical assumptions on its ingredients, this system can be seen to share with (\ref{CC1})
essential parts of the structural features discussed near (\ref{en}). 
In consequence, a spatially two-dimensional analysis can be designed similarly as that for (\ref{CC1}), leading to the
following result.
\begin{theo}\label{theo17} (\cite{jin_wang}) \quad
  Let $\Omega\subset\R^2$ be a bounded domain with smooth boundary, let $\chi>0$ and $\gamma>0$,
  and suppose that $f\in C^1([0,\infty)), g\in C^1([0,\infty))$ and $h\in C^2([0,\infty))$ are such that
  \be{17.1}
	f(0)=0,
	\quad f(s) \le -\theta s
	\mbox{ for all $s\ge 0$}
	\quad \mbox{and} \quad
	\frac{d}{ds} \big( \frac{f(s)}{s}\Big) \le - K_f
	\mbox{ for all $s> 0$,}
  \ee
  that
  \be{17.2}	
	g(0)=0,
	\quad
	g(s) >0
	\mbox{ for all $s\in (0,s_0)$}
	\quad \mbox{and} \quad
	g(s)<0
	\mbox{ for all $s>s_0$,}
  \ee
  that
  \be{17.3}
	h(0)=0
	\quad \mbox{as well as} \quad
	h(s)>0, \ h'(s)>0 \mbox{ and } h''(s) \le 0
	\mbox{ for all $s>0$,}
  \ee
  and that
  \be{17.4}
	\frac{g}{h}
	\mbox{ is continuously extensible to } [0,\infty).
  \ee
  with some $\theta>0$, $K_f>0, K_g>0$ and $s_0>0$.
  Then for arbitrary $u_0$ and $v_0$ fulfilling (\ref{init}), there exists a global classical solution of (\ref{CC6})
  which is bounded in the sense that with some $C>0$ we have
  \bas
	\|u(\cdot,t)\|_{L^\infty(\Om)} + \|v(\cdot,t)\|_{W^{1,\infty}(\Om)} \le C
	\qquad \mbox{for all } t>0.
  \eas
\end{theo}
Forming a typical representative of statements on large time behavior in multi-species parabolic systems involving
taxis, the following theorem asserts that the dynamics of the associated ODE system essentially determines the asymptotics
also for (\ref{CC6}) provided that either the former is essentially trivial, or the considered chemotactic coefficient
is appropriately small:
\begin{theo}\label{theo18} (\cite{jin_wang}) \quad
  Let $\Omega\subset\R^2$ be a bounded domain with smooth boundary, let $\gamma>0$,
  and assume that besides satisfying (\ref{17.1})-(\ref{17.4})
  with some $\theta>0$, $K_f>0, K_g>0$ and $s_0>0$, the functions $f, g$ and $h$ are such that
  \bas
	\lim_{s\searrow 0} \frac{g(s)}{h(s)} >0
	\qquad \mbox{and} \qquad
	\Big(\frac{g}{h}\Big)' <0
	\mbox{ on } (0,\infty).
  \eas
  i) \ If
  \bas
	\gamma h(s_0)<\theta,
  \eas
  then for any choice of $(u_0,v_0)$ fulfilling (\ref{init}), the solution $(u,v)$ of (\ref{CC6}) from Theorem \ref{theo17}
  satisfies
  \bas
	u(\cdot,t)\to 0
	\quad \mbox{and} \quad
	v(\cdot,t)\to s_0
	\quad \mbox{in } L^\infty(\Om)
	\qquad \mbox{as } t\to\infty.
  \eas
  ii) \ If
  \bas
	\gamma h(s_0)>\theta,
  \eas
  then there exists $\chi_0>0$ such that if $\chi\in (0,\chi_0)$, given any $(u_0,v_0)$ such that (\ref{init}) holds, 
  for the solution $(u,v)$ of (\ref{CC6}) found in Theorem \ref{theo17} we have
  \bas
	u(\cdot,t)\to u_\star
	\quad \mbox{and} \quad
	v(\cdot,t)\to v_\star
	\quad \mbox{in } L^\infty(\Om)
	\qquad \mbox{as } t\to\infty,
  \eas
  where $(u_\star,v_\star) \in (0,\infty) \times (0,s_0)$ is the unique solution of the algebraic system
  \bas
	\gamma u_\star h(v_\star) + f(u_\star)=0,
	\qquad
	- u_\star h(v_\star) + g(v_\star)=0.
  \eas
\end{theo}
Precedents and close relatives
concerned with questions of global solvability, boundedness and spatially homogeneous asymptotics
in (\ref{CC6}) and neighboring systems can be found in \cite{tao_prey}, \cite{sining_prey}, \cite{junping_shi},
\cite{tian_xiang}, \cite{JP_wang2}, \cite{jin_wang_yin}, \cite{dan_li}, \cite{li_yuxiang_preytaxis} and \cite{win_food_supported},
and an impression about the numerous recent developments in this field can be gained upon consulting
\cite{tello_wrzosek}, \cite{zhicheng_wang}, \cite{JP_wang}, \cite{ren1}, \cite{li_yuxiang_preytaxis},
\cite{inkyung_ahn}, \cite{wang_li_shi}, \cite{wrzosek2} and
\cite{liangchen_wang_prey}, for instance.\\[5mm]
From a mathematical point of view, considerably more thorough complexifications of (\ref{CC1}) are obtained when couplings
to further taxis mechanisms are included. As an exemplary case of such an extension, we here briefly consider the model
\be{CC7}
        \lball
	u_{1t} = \Delta u_1 - \nabla \cdot (u_1\nabla v),
	\qquad & x\in\Omega, \ t>0, \\[1mm]
	u_{2t} = \Delta u_2 - \nabla \cdot (u_2 \nabla u_1),
	\qquad & x\in\Omega, \ t>0, \\[1mm]
	v_t = \Delta v - (u_1+u_2)v - \mu v + r(x,t),
	\qquad & x\in\Omega, \ t>0, \\[1mm]
	\frac{\partial u_1}{\partial\nu}=\frac{\partial u_2}{\partial\nu}=\frac{\partial v}{\partial\nu}=0,
	\qquad & x\in\pO, \ t>0, \\[1mm]
	u_1(x,0)=u_{10}(x), \quad u_2(x,0)=u_{20}(x),  \quad v(x,0)=v_0(x),
	\qquad & x\in\Omega,
	\end{array} \right.
	\tag{CC8}
\ee
for the dynamics in so-called forager-exploiter systems.
Here, individuals in a first population, consisting of foragers and measured through the respective population density $u_1$,
are attracted by a food resource at concentration $v$, while a second population of exploiters, accordingly represented 
by $u_2$, orient their movement toward increasing concentrations of foragers.
Besides being consumed, nutrient is assumed to undergo spontaneous decay, and to be renewed from an external source
at a rate $r(x,t)$ (\cite{tania}).\abs
As the presence of the additional taxis process determining the evolution of $u_2$ apparently rules out
any global energy structure similar to that from (\ref{en}), already at the stage of questions from global solvability theory
an analysis of (\ref{CC7}) seems to require alternative tools. 
After all, by an appropriate further development of the approach outlined near (\ref{gradlogu}) it becomes possible
to derive a result on global generalized solvability in domains of arbitrary dimension, 
in a style similar to that from Theorem \ref{theo7}, 
under suitable smallness assumptions on both $v_0$ and $r$ (\cite{win_scrounge}; cf.~also \cite{li_yifu_wang}
for a related result on global existence even of classical small-data solutions).
Only in the case when $\Om\subset\R$ is a bounded interval, thanks to accordingly strong smoothing regularization features
of the Neumann heat semigroup, inter alia ensuring {\em a priori} bounds for $v_x$ with respect to the norm in $L^q(\Om)$
for arbitrary finite $q\ge 1$, global bounded solutions to (\ref{CC7}) have so far been constructed
for data of arbitrary size (\cite{taowin_scrounge}).\abs
Apart from that, the inclusion of additional relaxation mechanisms, such as saturation effects in nutrient consumption
or taxis,
nonlinear strengthening of diffusion at large densities of foragers and exploiters, or also logistic-type growth restrictions,
has been found to ensure global existence of solutions within various concepts of solvability also in higher-dimensional domains
(\cite{JP_wang_scrounge}, \cite{yuanyuan}, \cite{cao_scrounge}, \cite{black_scrounge}, \cite{JP_wang_scrounge2},
\cite{JP_wang_scrounge4}, \cite{chunlai_mu}, \cite{bin_liu_scrounge}, \cite{bin_liu_scrounge2}).
Solutions are known to stabilize to semi-trivial homogeneous equilibria, 
both in the one-dimensional version (\ref{CC7}) and in some of its modified variants
in higher-dimensional cases, if adequate assumptions are imposed either on smallness of the initial data or $r$, or on strength
of the respective relaxation mechanisms (\cite{cao_tao}, \cite{li_yifu_wang}, \cite{JP_wang_scrounge3},
\cite{zhongping_li}, \cite{zhongping_li2}, \cite{taowin_scrounge}). 
\mysection{Can chemotaxis-consumption models support spatial structures?}
In each of the results on large time behavior reported on above, the respective dynamics were determined
by convergence to homogeneous states. In the context of the models (\ref{CC1})-(\ref{CC5}), this may be viewed
as confirming the naive guess that in the long-term limit, all nutrient should be consumed and hence any taxis gradient
asymptotically disappear; in fact, also in the complex settings of (\ref{CC6}) and (\ref{CC7}) all the mentioned 
result on stabilization have been derived under assumptions leading to essentially similar situations.\abs
This final section attempts to describe some developments which at their core aim at indicating that despite
this homogenizing role played by consumption, nevertheless some emergence and stabilization of spatial structures 
can be supported.
\subsection{Effects of nonhomogeneous boundary conditions}
Our first focus in this regard will be on possible influences of boundaries; indeed, already in the experimental 
framework underlying the introduction of the simple model (\ref{CC1}), some nontrivial effects related
to nutrient balancing at boundaries have been pointed out (\cite{goldstein2004}, \cite{goldstein2005}).\abs
In fact, a corresponding consideration on refined modeling near boundaries (\cite{braukhoff17,braukhoff_lankeit,braukhoff_tang}) suggests
the variant of (\ref{CC1}) given by
\be{CC8}
	\lball
	u_t=\Del u - \na\cdot (u\na v), 
	\qquad & x\in\Om, \ t>0, \\[1mm]
	v_t = \Del v - uv
	\qquad & x\in\Om, \ t>0, \\[1mm]
	\frac{\pa u}{\pa\nu}=0, \quad \frac{\pa v}{\pa\nu}=\beta(x) (\gamma(x)-v),
	\qquad & x\in\pO, \ t>0, \\[1mm]
	u(x,0)=u_0(x), \quad v(x,0)=v_0(x),
	\qquad & x\in\Om,
	\ear
	\tag{CC9}
\ee
as a modification which takes into account these boundary effects in a realistic manner; here,
$\beta$ and $\gamma$ are given nonnegative functions on $\bom$. In place of these Robin-type 
boundary conditions, inhomogeneous Dirichlet conditions, $v|_{\pO}=\gamma(x)$, may be imposed (e.g. \cite{goldstein2005,peng_xiang_19}),
which arise as limiting case of the former, cf. \cite[Prop.~5.1]{braukhoff_lankeit} \abs
A treatment along the lines of \eqref{en}--\eqref{CC1.1} meets the difficulty that both in \eqref{D} and in \eqref{FE} 
boundary integrals appear which no longer can be immediately estimated.

As an observation of major importance for the development of a basic solvability, the authors in \cite{braukhoff_tang} 
incorporated a boundary contribution into the functional and 
derived a quasi-energy inequality of the form
\bas
	\F_4'(t) + a\F_4(t) \le b
	\qquad \mbox{for all } t>0,
\eas
where with some $L>0$,		
\bea{en4}
	\F_4(t)
	&:=& \io u\ln u + \frac{1}{2} \io \frac{|\na v|^2}{v}
	+ \int_{\pO} \beta(x) \cdot \Big\{ \gamma(x) \ln \frac{\gamma(x)}{v} - \gamma(x)+v\Big\} \nn\\
	& & + L \io \Big\{ v\ln \frac{v}{\gamma(x)} - v + \gamma(x)\Big\},
	\qquad t\ge 0.
\eea
Accordingly obtained {\em a priori} estimates led to the following analogue of the solvability statement in Theorem \ref{theo1}.
\begin{theo}\label{theo19} (\cite{braukhoff_tang}) \quad
  Let $\Om\subset\R^2$ be a bounded domain with smooth boundary, let $\beta\in C^1(\bom)$ and $\gamma\in C^1(\bom)$
  be positive, and suppose that (\ref{init}) holds.
  Then there exists a global classical solution of (\ref{CC8}), uniquely determined through (\ref{reg}),
  such that $u\ge 0$ and $v>0$ in $\bom\times (0,\infty)$, and that $\sup_{t>0} \|u(\cdot,t)\|_{L^\infty(\Om)} <\infty$.
\end{theo}
We remark here that in three- and higher-dimensional domains, a corresponding counterpart of Theorem \ref{theo2}
can be derived on the basis (\ref{en3}) and a straightforward adaptation of the argument from Theorem \ref{theo2};
for the particular case $n=3$, this can be found detailed in \cite{braukhoff_tang}.
A result on global solvability in a parabolic-elliptic of (\ref{CC8}) involving tensor-valued sensitivities has recently been achieved
in \cite{ahn_kang_lee} by means of an interesting method suitably decomposing the spatial domain.\abs
In \cite{WWX_CPDE}, a localized energy functional (cutting off the functions near the boundary) served as key ingredient in the construction of global generalized 
solutions to the Dirichlet variant of the system, whereas \cite{lankeit_win_NON} relied on the de facto one-dimensionality of said problem in a radially symmetric setting in order 
to derive sufficient bounds on the integrands in the boundary terms. 
\abs
Now with respect to the large time behavior, it is quite evident from (\ref{CC8}) 
that for any choice of positive $\beta$ and $\gamma$, no solution can exhibit a behavior similar to that in (\ref{conv});
indeed, any stabilization toward a limit function $v_\infty$ in the second component should require $v_\infty$ to be nonconstant.
In the case when $\gamma$ is constant, a more thorough description of the set of equilibria can be given as follows.
\begin{theo}\label{theo20} (\cite{braukhoff_lankeit}) \quad
  Let $n\ge 1$ and $\Om\subset\R^n$ be a bounded domain with smooth boundary, and 
  let $\beta\in \bigcup_{\theta>0} C^{1+\theta}(\bom)$ 
  and $\gamma\in (0,\infty)$.
  Then for any $M>0$ there is exactly one pair $(u,v)\in (C^1(\bom)\cap C^2(\Om))^2$ such that $\io u=M$, and that
  $\bom\times [0,\infty)\ni (x,t)\mapsto (u(x),v(x))$ forms a classical solution of the boundary-value problem in (\ref{CC8}).
  This solution satisfies $u>0$ and $v>0$ in $\bom$ as well as 
  \bas
	u\not\equiv const.
	\qquad \mbox{and} \qquad
	v\not\equiv const.
  \eas
\end{theo}
Corresponding results for \eqref{CC8} with Dirichlet condition can be found in \cite{lankeit_win_NON} and \cite{lee_wang_yang}, 
the latter of which is concerned with the formation of boundary layers in the limit of vanishing diffusion of $v$.\abs
An interesting open problem now consists in establishing a link between solutions of the parabolic problem (\ref{CC8}) 
and the steady states just characterized; in view of the uniqueness statement in Theorem \ref{theo20} it may be conjectured
that $\omega$-limit sets at least of smooth trajectories should consist exclusively of such equilibria at the respective mass level.
A rigorous confirmation of this has recently been achieved for certain small-data solutions to a parabolic-elliptic 
variant of (\ref{CC8}) (\cite{fuest_lankeit_mizukami}).\abs
For some further system variants closely related to \eqref{CC8}, results on global solvability were obtained in \cite{knosalla1}, \cite{knosalla3}, 
\cite{black_winkler_22}, \cite{tian_xiang}, \cite{wu_xiang_22}, \cite{black_wu1}, \cite{black_wu2},  
and corresponding steady-state analysis can be found in \cite{knosalla2}, \cite{knosalla4}.
\abs 
Two recent examples of quasilinear variants of (\ref{CC8}) involving nonhomogeneous Dirichlet boundary conditions and an elliptic equation for $v$ 
show that the combination of signal consumption with repulsive chemotaxis may even lead to finite-time blow-up
(\cite{wang_winkler}, \cite{ahn_winkler}). 
\subsection{The role of cross-degeneracies in structure formation}
In situations when significant parts of structure formation are expected to occur in the interior of a spatial domain,
model modifications merely concentrating on alternative boundary conditions may be insufficient;
striking experimental findings on fractal-like patterning in populations of {\em Bacillus subtilis} suspended to 
nutrient-poor agars form particularly simple but convincing examples (\cite{fujikawa}, \cite{fujikawa_matsushita}).
A recent modeling approach proposes systems of the form
\be{CC9}
	\lball
	u_t=\na\cdot (uv\na u) - \na \cdot (u^2 v\na v),
	\qquad & x\in\Om, \ t>0, \\[1mm]
	v_t = \Del v - uv
	\qquad & x\in\Om, \ t>0, \\[1mm]
	(uv\na u - \chi u^2 v\na v)\cdot \nu=\na v\cdot\nu=0,
	\qquad & x\in\pO, \ t>0, \\[1mm]
	u(x,0)=u_0(x), \quad v(x,0)=v_0(x),
	\qquad & x\in\Om,
	\ear
	\tag{CC10}
\ee
for the description of such phenomena (\cite{plaza}), where a key novelty consists in the introduction of the factor $v$ 
not only to the cross-diffusion term, but especially also to the part related to random diffusion of cells.
The reduction of motility thereby accounted for seems in good accordance with experimentally gained knowledge
on bacterial migration in nutrient-poor environments; 
however, the consequences for the mathematical properties of the resulting evolution system appear to be rather drastic.
In fact, the cross-degeneracy of diffusion thereby induced does not only go along with an apparent
loss of favorable parabolic smoothing features, but it seems to bring about new qualitative features 
which let spatial structures seem to become an essential core, rather than an exception:
\begin{theo}\label{theo21} (\cite{liwin1}, \cite{win_TRAN}) \quad
  Let $\Om\subset\R$ be a bounded open interval, and assume that $u_0\in W^{1,\infty}(\Om)$ is nonnegative with
  $u_0\not\equiv 0$.\abs
  i) \ Whenever $v_0\in W^{1,\infty}(\Om)$ is positive in $\bom$, there exist functions
  \bas
	\lbal
	u\in C^0(\bom\times [0,\infty)) \cap L^\infty(\Om\times (0,\infty)) \qquad \mbox{and} \\[1mm]
	v\in C^0(\bom\times [0,\infty)) \cap C^{2,1}(\bom\times (0,\infty)),
	\ear
  \eas
  such that
  \bas
	u^2 v \mbox{ and } u^2 v_x \mbox{ belong to } \in L^1_{loc}(\bom\times [0,\infty)),
  \eas
  that $u\ge 0$ and $v>0$ in $\bom\times (0,\infty)$, and that (\ref{CC9}) is solved in the sense that (\ref{2.3})
  holds for all $\vp\in C_0^\infty(\bom\times [0,\infty))$, and that
  \bas
	-\int_0^\infty \io u\vp_t
	- \io u_0 \vp(\cdot,0)
	= \frac{1}{2} \int_0^\infty \io u^2 v_x \vp_x
	+ \frac{1}{2} \int_0^\infty \io u^2 v \vp_{xx} 
	+ \int_0^\infty \io u^2 v v_x \vp_x
  \eas
  for all $\vp\in C_0^\infty(\bom\times [0,\infty))$ fulfilling $\vp_x=0$ on $\pO\times (0,\infty)$.
  Moreover, there exists $u_\infty\in C^0(\bom)$ such that
  \be{21.1}
	u(\cdot,t) \to u_\infty
	\ \mbox{in } L^\infty(\Om)
	\quad \mbox{and} \quad
	v(\cdot,t) \to 0
	\ \mbox{in } L^\infty(\Om)
	\qquad \mbox{as } t\to\infty.
  \ee
  ii) \ There exist $\theta>0$ and $\del>0$ such that if $0<v_0\in W^{1,\infty}(\Om)$ with
  \be{21.2}
	\|v_0^\theta\|_{W^{1,\infty}(\Om)} \le \del,
  \ee
  then in (\ref{21.1}) we have 
  \be{21.3}
	u_\infty \not\equiv const.
  \ee
\end{theo}
We remark that beyond this result on potentially nontrivial stabilization, the cross-degeneracy in (\ref{CC9}) also 
facilitates a stability property that apparently is not frequently seen in contexts of parabolic problems.
Namely, it is evident that for arbitrary suitably regular $u_\star:\Om\to [0,\infty)$, the pair $(u_\star,0)$
formally constitutes a steady-state solution of (\ref{CC9}); 
it can now be shown that within suitable topologies, {\em any} of these equilibria is stable in the classical flavor:
Trajectories emanating from initial data nearby remain close throughout evolution (\cite{win_TRAN}).
In this sense, the cross-degeneracy in diffusion does not only bring about nontrivial asymptotics as in (\ref{21.1})-(\ref{21.3}),
but also entails stabilization of arbitrary structures in nutrient-poor settings for which (\ref{CC9}) has been designed.\abs
In higher-dimensional domains, up to now only certain small-signal solutions of (\ref{CC9})
have been shown to exist globally (\cite{win_NONRWA}); for initial data of arbitrary size,
two results on global solvability in two- and three-dimensional settings have recently been derived for
a variant of (\ref{CC9}) involving some appropriate saturation of taxis at large population densities (\cite{genglin_JDE}).\abs
Similar effects of such cross-degeneracies could recently be detected to occur also in the model
\be{CC10}
	\lball
	u_t=\Del \big( u\phi(v) \big),
	\qquad & x\in\Om, \ t>0, \\[1mm]
	v_t = \Del v - uv
	\qquad & x\in\Om, \ t>0, \\[1mm]
	\na \big( u\phi(v) \big)\cdot\nu=\na v\cdot\nu=0,
	\qquad & x\in\pO, \ t>0, \\[1mm]
	u(x,0)=u_0(x), \quad v(x,0)=v_0(x),
	\qquad & x\in\Om,
	\ear
	\tag{CC11}
\ee
for microbial motion influenced by certain local sensing mechanisms (\cite{cho_kim}, \cite{fu}).
Indeed, while solutions behave essentially as described in (\ref{conv}) when here the key ingredient
$\phi$ is sufficiently regular and strictly positive on $[0,\infty)$ (\cite{laurencot}, \cite{li_winkler}; see also  \cite{li_zhao}, \cite{yoon}
for precedent existence results, and 
\cite{wenbin11}, \cite{li_yan}, \cite{liangchen_wang2} as well as \cite{li_wang_pan,qiu_mu_tu} for findings on some close relatives),
a statement similar to those in (\ref{21.1})-(\ref{21.3}) holds when $\phi$ gives rise to a cross-degeneracy by
satisfying $\phi(s) \simeq s$ near $s=0$ in a suitable sense (\cite{win_ANIHPC}).\abs
Another local sensing mechanism forms the basis of the system
\be{CC12}
	\lball
	u_t=\nabla\cdot \big( \frac{1}v \nabla \big(\frac{u}{v}\big) \big),
	\qquad & x\in\Om, \ t>0, \\[1mm]
	v_t = \eta\Del v - uv,
	\qquad & x\in\Om, \ t>0, \\[1mm]
	\na \big( u\phi(v) \big)\cdot\nu=\na v\cdot\nu=0,
	\qquad & x\in\pO, \ t>0, \\[1mm]
	u(x,0)=u_0(x), \quad v(x,0)=v_0(x),
	\qquad & x\in\Om,
	\ear
	\tag{CC12}
\ee
which stems from a random walk where steps have constant length if distance is measured with respect to the 'food-metric' (\cite{choi_kim}).  
For \eqref{CC12} (with $\eta=0$ or $\eta\ge0$) in one-dimensional domains, traveling waves (\cite{choi_kim}), global existence (\cite{ahn_choi_yoo}, \cite{ahn_choi_yoo2})  and first results indicating structure formation were found \cite{song_li}.

\bigskip
{\bf Acknowledgement.} \quad
  The second author acknowledges support of the {\em Deutsche Forschungsgemeinschaft} (Project No.~462888149).

{\footnotesize 
  \setlength{\parskip}{0pt}
  \setlength{\itemsep}{0pt plus 0.05ex}
\setlength{\baselineskip}{0.92\baselineskip}

 }
 \end{document}